\documentstyle[12pt,amsfonts,amssymb,leqno]{article}

\newtheorem{proposition}{Proposition}[section]
\newtheorem{thm}[proposition]{Theorem}
\newtheorem{cor}[proposition]{Corollary}
\newtheorem{lemma}[proposition]{Lemma}
\newtheorem{defn}[proposition]{Definition}

\newcommand{\forces}{\ \ |\!\!\!|\!- }

\begin{document}
\setcounter{section}{-1}

\begin{center}
{\Large \bf $L({\mathbb R})$ absoluteness under proper forcings}
\end{center}

\begin{center}
\renewcommand{\thefootnote}{\fnsymbol{footnote}}
\renewcommand{\thefootnote}{arabic{footnote}}
\renewcommand{\thefootnote}{\fnsymbol{footnote}}
{\large Ralf-Dieter Schindler}${}^{a}$\footnote[1]{The author is indebted to Juan
Bagaria, Sy Friedman, and Philip Welch for stimulating hints and observations.

\noindent
1991 {\it Mathematics Subject Classification.} Primary 03E55, 03E15. Secondary
03E35, 03E60.\\
Keywords: set theory/descriptive set theory/proper forcing/large cardinals.}
\renewcommand{\thefootnote}{arabic{footnote}}
\end{center}
\begin{center} 
{\footnotesize
${}^a${\it Institut f\"ur Logistik, Universit\"at Wien, 1090 Wien, Austria}} 
\end{center}

\begin{center}
{\tt rds@logic.univie.ac.at}\\
\end{center}

\begin{abstract}
\noindent
We isolate a new large cardinal concept, "remarkability." It turns out that
the existence of a remarkable cardinal is 
equiconsistent with $L({\mathbb R})$ absoluteness under proper forcings. As
remarkable cardinals are
compatible with $V = L$, this means that said absoluteness does not imply $\Pi^1_1$
determinacy.
\end{abstract}

\section{Introduction.}

It is well-known that the existence of large cardinals implies various forms of
determinacy. It is also true that the existence of large cardinals
implies that the theory of $L({\mathbb R})$ cannot be changed by set-forcing.
It was an interesting observation (due to Steel and Woodin) to see
that already the fact that
the theory of $L({\mathbb R})$ cannot be changed by set-forcing
implies that the axiom of determinacy holds in $L({\mathbb R})$.

Large cardinals imply more: \cite{NZ1} and \cite{NZ2} show that
-- under appropriate assumptions -- the "boldface" theory of 
$L({\mathbb R})$ cannot be changed by set-sized {\em proper} forcing. Here, "boldface"
means that reals from the ground model as well as ordinals are allowed as parameters. 
A natural question arises: $(\star)$
which amount of determinacy do the conclusions of the main
theorems of \cite{NZ1} and \cite{NZ2} give (back)?

This question is particularily interseting, as the forcing which Steel and Woodin use
to prove their above-mentioned observation collapses $\omega_1$. The question thus
really is whether the Steel-Woodin argument can be refined by using more "coding like"
forcings instead. \cite{CKRF} gave some partial answers, albeit in a somewhat
different direction. We here provide a straight answer to $(\star)$.
In order to formulate it let us introduce some terminology.

\begin{defn}\label{def-absoluteness}
Let ${\cal F} \subset V$ be a class of posets. We say that $L({\mathbb R})$ {\em is
absolute under forcings of type} ${\cal F}$ if for all posets $P \in {\cal F}$, for all
$H$ being $P$-generic over $V$, for all formulae $\Phi({\vec v})$, and for all ${\vec x}
\in {\mathbb R}^V$ do we have that
$$L({\mathbb R}^V) \models \Phi({\vec x}) \Leftrightarrow 
L({\mathbb R}^{V[G]}) \models \Phi({\vec x}).$$
We say that $L({\mathbb R})$ {\em is
absolute under proper forcings} if $L({\mathbb R})$ is
absolute under forcings of type ${\cal F}$ where ${\cal F} = \{ P \in V \colon P$ is
proper $\}$.
\end{defn}

\begin{defn}\label{def-embedding-thm}
Let ${\cal F} \subset V$ be a class of posets. We say that {\em the}
$L({\mathbb R})$ {\em embedding theorem holds for
forcings of type} ${\cal F}$ if for all posets $P \in {\cal F}$, for all
$H$ being $P$-generic over $V$, for all formulae $\Phi({\vec v})$, for all ${\vec
\alpha} \in OR$, and for all ${\vec x}
\in {\mathbb R}^V$ do we have that
$$L({\mathbb R}^V) \models \Phi({\vec \alpha},{\vec x}) \Leftrightarrow 
L({\mathbb R}^{V[G]}) \models \Phi({\vec \alpha},{\vec x}).$$
We say that {\em the} $L({\mathbb R})$ {\em embedding theorem holds for
proper forcings} if $L({\mathbb R})$ is
absolute under forcings of type ${\cal F}$ where ${\cal F} = \{ P \in V \colon P$ is
proper $\}$.
\end{defn}

\begin{defn}\label{def-anti-coding}
Let ${\cal F} \subset V$ be a class of posets. We say that {\em the}
$L({\mathbb R})$ {\em anti coding theorem holds for
forcings of type} ${\cal F}$ if for all posets $P \in {\cal F}$, for all
$H$ being $P$-generic over $V$, and for all $A \subset OR$ with $A \in V$
do we have that
$$A \in L({\mathbb R}^V) \Leftrightarrow 
A \in L({\mathbb R}^{V[G]}).$$
We say that {\em the} $L({\mathbb R})$ {\em anti coding theorem holds for
proper forcings} if the $L({\mathbb R})$ anti coding theorem holds
for forcings of type ${\cal F}$ where ${\cal F} = \{ P \in V \colon P$ is
proper $\}$.
\end{defn}

Our main theorem, \ref{equicon}, will then say that $L({\mathbb R})$ absoluteness
under proper forcings is equiconsistent with the existence of
what we call a "remarkable" cardinal, and the same holds for
the $L({\mathbb R})$ embedding theorem for proper forcings,
as well as the $L({\mathbb R})$ anti coding
theorem for proper forcings. 
As remarkable cardinals turn out to be compatible with $V=L$, this means that the
answer to $(\star)$ is that 
the conclusions of the main
theorems of \cite{NZ1} and \cite{NZ2} do not even imply $\Pi^1_1$ determinacy.

We hope that there will be other applications of remarkable cardinals in the future;
anyway, "remarkability" seems to be an interesting concept.

\section{Preliminaries.}

\begin{lemma}\label{ex-of-emb1}
Let ${\cal M} = (M;(R_i \colon i<n))$ and ${\cal N} = (N;(S_i \colon i<m))$ be models
such that $n \leq m$, $R_i$ has the same arity as $S_i$ for $i<n$, and $M$ is countable.
Then there is a tree $T$ 
of height $\leq \omega$ searching for $(R_i \colon n \leq i < m)$
together with an
elementary embedding $$\pi \colon
(M;(R_i \colon i<m)) \rightarrow (N;(S_i \colon i<m)).$$
\end{lemma}

{\sc Proof.} Let $(e_i \colon i<\omega)$ be an enumeration of $M$, and let 
$(\Phi_i({\vec
v}) \colon i<\omega)$ be an enumeration of all formulae of the language associated with
${\cal N}$. Let $\sharp(i)$ denote the arity of $R_i$ ($=$ of $S_i$) for $i<n$.
Let $\gamma \colon \omega \rightarrow \colon \omega \times {}^{<\omega}\omega$ 
be such that $\Phi_{\gamma(i)_0}$ has the variables with indices $< dom(\gamma(i)_1)$
as its free variables and $ran(\gamma(i)_1) \subset i-1$, 
and such that
$\gamma$ is "onto" in the obvious sense. Let ${\cal F}$
be a Skolem function for ${\cal N}$; more precisely, let ${\cal F}(i,{\vec x})$ be
such that $${\cal N} \models \exists y \ \Phi_i(y,{\vec x}) \Rightarrow
\Phi_i({\cal F}(i,{\vec x}),{\vec x})$$ (if there is no such $y$ then
we let ${\cal F}(i,{\vec x})$ undefined). 
Let the $k^{th}$ level
of $T$ consist of sequences $f \colon k \rightarrow N$ 
such that $f \upharpoonright k-1 \in (k-1)^{st}$ level of $T$, 
$$\forall i<n \forall \{
l_1,...,l_{\sharp(i)} \} \subset k \ ( 
\ R_i(e_{l_1},...,e_{l_{\sharp(i)}}) \Leftrightarrow 
S_i(f(l_1),...,f(l_{\sharp(i)})) \ ) {\rm , \ and }$$
$$f(k-1) = {\cal F}(\gamma(k)_0,f \circ \gamma(k)_1(1),..., f \circ
\gamma(k)_0,\gamma(k)_1(dom(\gamma(k)_1)-1))$$ (if this is defined, otherwise
we let $f(k-1) =$ an arbitrary element of $N$).

Now if $f \colon \omega \rightarrow N$ is given by an infinite branch through $T$ then
it is easy to see that setting $R_i(e_{l_1},...,e_{l_p}) \Leftrightarrow 
S_i(f(l_1),...,f(l_p))$ for $n \leq i <m$ and $\pi(e_i) = f(i)$ we get relations and
an embedding as desired. On the other hand, any such relations together with some such
embedding defines an infinite branch through $T$.

\bigskip
\hfill $\square$ (\ref{ex-of-emb1})

\bigskip
As an immediate corollary to this proof we get the following.

\begin{lemma}\label{ex-of-emb2}
Let ${\cal M} = (M;(R_i \colon i<n))$ and ${\cal N} = (N;(S_i \colon i<m))$ be models
such that $n \leq m$, $R_i$ has the same arity as $S_i$ for $i<n$, and $M$ is countable.
Let $Q$ be an admissible set such that ${\cal M}$, ${\cal N} \in Q$, and 
$M$ is countable in $Q$. If in $V$ there are $R_i$, $n \leq i < m$,
together with an
elementary embedding $$\pi \colon (M;(R_i \colon i<m)) \rightarrow 
(N;(S_i \colon i<m))$$
then such $R_i$, $\pi$ also exist in $Q$.
\end{lemma}
\section{Remarkable cardinals.}

\begin{defn}\label{def-remarkable}
A cardinal $\kappa$ is called {\em remarkable} iff for all regular cardinals $\theta >
\kappa$ there is some $\pi \colon M \rightarrow H_\theta$ with $M$ being
countable and transitive and $\kappa
\in ran(\pi)$ and 
such that, setting ${\bar \kappa} = \pi^{-1}(\kappa)$,
there is some $\sigma \colon M \rightarrow N$ with 
critical point ${\bar \kappa}$ and such that
$N$ is
countable and transitive, ${\bar \theta} = M \cap OR$ is a regular cardinal in $N$, 
$M = H^N_{\bar \theta}$, and $\sigma({\bar \kappa}) > {\bar \theta}$. 
\end{defn}

As a matter of fact, "remarkability" relativizes down to $L$, i.e.,
any remarkable cardinal is also remarkable in
$L$ (cf. \ref{rel-down-to-L} below). 
Hence the existence of remarkable cardinals is consistent with $V = L$.
One can also show that every remarkable cardinal is totally indescribable. 
In particular, the least measurable cardinal is not remarkable. However, every strong
cardinal {\em is} remarkable, and we shall see below (cf. \ref{ind-rem})
that every Silver indiscernible
is remarkable in $L$. 

The following two lemmata will give information as to where remarkable cardinals
sit in the large cardinal hierarchy.  

\begin{lemma}\label{erdos-rem}
Let $\kappa \rightarrow (\omega)^{<\omega}$. Then 
there are $\alpha < \beta < \omega_1$ such that
$L_\beta \models$ "$ZFC \ + \ \alpha$ is a
remarkable cardinal."
\end{lemma}

{\sc Proof.} We may assume that $V=L$, as $\kappa \rightarrow (\omega)^{<\omega}$
relativizes down to $L$. 
Let $\pi \colon L_\gamma \rightarrow L_\kappa$ be 
such that $ran(\pi)$ is
the Skolem hull in $L_\kappa$ of $\omega$ many indiscernibles for $L_\kappa$.
Let $\alpha$, $\beta$ (with $\alpha < \beta$) 
be the images of the first two indiscernibles under $\pi^{-1}$.
Of course, $L_\beta \models ZFC$, as any of the indiscernibles in inaccessible in $L$.
We claim that $\alpha$ is remarkable in $L_\beta$.

Let $\theta < \beta$ be regular in $L_\beta$ with $\theta > \alpha$. 
There is $\sigma \colon L_\gamma
\rightarrow L_\gamma$ with $\sigma(\alpha) = \beta$, obtained from shifting the
indiscernibles. I.e., there is some countable $L_{\bar \theta}$ (namely,
$L_\theta$) together with some ${\bar \pi} \colon L_{\bar \theta} \rightarrow 
L_{\pi(\theta)}$ (namely, $\pi \upharpoonright L_\theta$)
such that $\pi(\alpha)$ is in the range of ${\bar \pi}$ and
there is some ${\bar 
\sigma} \colon L_{\bar \theta} \rightarrow 
L_{\tilde
\theta}$ (namely, $\sigma \upharpoonright L_\theta$)
with critical point ${\bar
\pi}^{-1}(\pi(\alpha))$
such that ${\tilde
\theta}$ is countable, ${\bar \theta}$ is a regular cardinal in $L_{\tilde
\theta}$, and ${\bar \sigma}({\bar
\pi}^{-1}(\pi(\alpha))) > {\bar \theta}$. 
As $\pi(\beta)$ is inaccessible in $L$, 
the same holds in $L_{\pi(\beta)}$. Pulling it back via 
$\pi^{-1}$ we
get that in $L_\beta$ there is some countable $L_{\bar \theta}$ 
together with some ${\bar \pi} \colon L_{\bar \theta} \rightarrow 
L_\theta$ such that $\alpha$ is in the range of ${\bar \pi}$ and
there is some ${\bar 
\sigma} \colon L_{\bar \theta} \rightarrow 
L_{\tilde
\theta}$ with critical point ${\bar
\pi}^{-1}(\alpha)$
such that ${\tilde
\theta}$ is countable, ${\bar \theta}$ is a regular cardinal in $L_{\tilde
\theta}$, and ${\bar \sigma}({\bar
\pi}^{-1}(\alpha)) > {\bar \theta}$. As $\theta > \alpha$ was an arbitrary regular
cardinal in $L_\beta$, we have shown that $\alpha$ is remarkable in $L_\beta$.

\bigskip
\hfill $\square$ (\ref{erdos-rem})

\bigskip
As an immediate corollary to this proof we get:

\begin{lemma}\label{ind-rem}
Suppose that $0^\sharp$ exists. Then every
Silver indiscernible is remarkable in $L$.
\end{lemma}

{\sc Proof.} 
A slight variation of the previous proof gives that $L_\beta \models$ "$\alpha$ is 
remarkable" whenever $\alpha < \beta$ are both indiscernibles for $L$. But then 
every
Silver indiscernible is remarkable in $L$.

\bigskip
\hfill $\square$ (\ref{ind-rem})

\begin{lemma}\label{rem-ineffable}
Let $\kappa$ be remarkable. Then 
there are $\alpha < \beta < \omega_1$ such that
$L_\beta \models$ "$ZFC \ + \ \alpha$ is a
ineffable cardinal."
\end{lemma}

{\sc Proof.} Let $\theta = \kappa^+$, and let $\pi$, $M$, $\sigma$, and $N$ be as in
\ref{def-remarkable}. Let $\alpha = \pi^{-1}(\kappa)$ and let $\beta =
\sigma(\alpha)$. It is easy to see that $L_\beta \models ZFC$. We claim that $\alpha$
is ineffable in $L_\beta$. 

Let $(A_i \colon i<\alpha) \in L_\beta$ be such that $A_i \subset i$ for all
$i<\alpha$, and let $C \in L_\beta$ be club in $\alpha$. There is $(A_i \colon \alpha
\leq i < \beta)$ such that $\sigma((A_i \colon i<\alpha)) = (A_i \colon i<\beta)$.
Notice that $A_\alpha \in N$, as ${\cal P}(\alpha) \cap N = {\cal P}(\alpha) \cap M$
by the properties of $M$, $\sigma$, and $N$. Now of course
$A_\alpha = \sigma(A_\alpha) \cap \alpha$, and also $\alpha \in \sigma(C)$.
This gives that 
$\alpha \in \sigma(\{ i<\alpha \colon A_i = A_\alpha \cap i \}) 
\cap \sigma(C)$, and thus via $\sigma$ we have that
$\{ i<\alpha \colon A_i = A_\alpha \cap i \} \cap C \not= \emptyset$.
As $C$ was arbitrary, we have shown that $\alpha$
is ineffable in $L_\beta$. 

\bigskip
\hfill $\square$ (\ref{rem-ineffable})

\bigskip
We now turn towards a useful characterization of remarkability.

\begin{defn}\label{def-char}
Let $\kappa$ be a cardinal. Let $G$ be
$Col(\omega,<\kappa)$-generic over $V$, let $\theta > \kappa$ be a
regular cardinal, and let
$X \in [H_\theta^{V[G]}]^\omega$. We say that $X$ {\em condenses remarkably} if 
$X = ran(\pi)$ for some 
elementary $$\pi \colon (H_\beta^{V[G \cap
H_\alpha^V]};\in,H_\beta^V,G \cap
H_\alpha^V) \rightarrow (H_\theta^{V[G]};\in,H_\theta^V,G)$$ where 
$\alpha = crit(\pi) < \beta < \kappa$ and $\beta$ is a regular cardinal.
\end{defn}

Notice that in the situation of \ref{def-char} we will have that $\alpha$ is
inaccessible in $V$, $G \cap H_\alpha^V$ is $Col(\omega,<\alpha)$-generic over $V$,
and hence $\beta$ is a regular cardinal in $V[G \cap H_\alpha^V]$, too.

\begin{lemma}\label{char-remarkable}
A cardinal $\kappa$ is remarkable if and only if for all regular cardinals $\theta >
\kappa$ do
we have that
$$\forces_{Col(\omega,<\kappa)}^V \ {\rm " } 
\{ X \in [H_{\check \theta}^{V[{\dot G}]}]^\omega \colon
X {\rm \ condenses \ remarkably } \} {\rm \ is \ stationary." }$$
\end{lemma}

{\sc Proof.} "$\Rightarrow$." 
Let $\kappa$ be remarkable, and let $\theta > \kappa$ be a regular cardinal. 
We may pick $\pi \colon M \rightarrow
H_{\theta^{+}}$ as in \ref{def-remarkable}, but with $\theta^+$ playing the role of
$\theta$. Let ${\bar \kappa}$, ${\bar \theta} = \pi^{-1}(\kappa$,
$\theta)$, 
and let $\sigma \colon M \rightarrow N$ with critical point ${\bar \kappa}$
be such that
$N$ is countable and transitive, $\rho = M \cap OR$ is regular in $N$,
$M = H^N_\rho$, and $\sigma({\bar \kappa}) > \rho$.
In $V$, we may pick $G$ being $Col(\omega,<{\bar \kappa})$-generic over $M$ (and hence
over $N$), and we may pick $G' \supset G$ being
$Col(\omega,<\sigma({\bar \kappa}))$-generic over $N$. We then
have that $\sigma$ naturally extends to ${\tilde \sigma} 
\colon M[G] \rightarrow N[G']$.

Let ${\cal M} = (H^{M[G]}_{\bar \theta};\in,H^M_{\bar \theta},G,(R_i \colon i<n) \in 
M[G]$ be any model of finite type.
Notice that ${\cal M} \in N[G']$ and is countable
there. 
By the existence of ${\tilde \sigma} \upharpoonright H^{M[G]}_{\bar \theta} \colon
{\cal M} \rightarrow {\tilde \sigma}({\cal M})$ together with \ref{ex-of-emb2}, 
we get that 
in $N[G']$ there is an elementary embedding $\tau$ of ${\cal M}$ into 
${\tilde \sigma}({\cal M})$.
This means that $$N[G'] \models 
\exists \alpha < \beta < \sigma({\bar \kappa}) 
\exists \tau \ ( \tau \colon (H_\beta^{V[G' \cap
H_\alpha]};\in,...) \rightarrow {\tilde \sigma}({\cal M}) 
\wedge \beta {\rm \ is \ regular \ }
).$$ Pulling this back via ${\tilde \sigma}$ gives that  
$$M[G] \models 
\exists \alpha < \beta < {\bar \kappa} 
\exists \tau \ ( \tau \colon (H_\beta^{V[G \cap
H_\alpha]};\in,...) \rightarrow {\cal M} \wedge \beta {\rm \ is \ regular \ }
).$$
As ${\cal M}$ was arbitrary, 
we have shown that $$\forces^M_{Col(\omega,<{\bar \kappa})} \
\{ X \in [H^{M[{\dot G}]}_{\check {\bar \theta}}]^\omega 
\colon X {\rm \ condenses \ remarkably \
} \} {\rm \ is \ stationary." }$$
Lifting this up via $\pi$ gives  
$$\forces_{Col(\omega,<\kappa)}^V \
\{ X \in [H_{\check \theta}^{V[{\dot G}]}]^\omega 
\colon X {\rm \ condenses \ remarkably \
} \} {\rm \ is \ stationary." }$$
As $\theta$ was arbitrary, this proves "$\Rightarrow$."

"$\Leftarrow$." Let $\theta > \kappa$ be a regular cardinal, and suppose that
$$\forces_{Col(\omega,<\kappa)}^V \
\{ X \in [H_{{\check \theta}}^{V[{\dot G}]}]^\omega 
\colon X {\rm \ condenses \ remarkably \
} \} {\rm \ is \ stationary." }$$
Let ${\bar \pi} \colon {\bar M} \rightarrow H_{\theta^{+}}$ with $M$ 
countable and transitive be such that $\kappa$, $\theta \in ran({\bar \pi})$.
Let ${\bar \kappa}$, ${\bar \theta} = {\bar \pi}^{-1}(\kappa$, $\theta)$.
In $V$, we may pick $G$ being $Col(\omega,<{\bar \kappa})$-generic over ${\bar M}$.
Because $$\forces^{\bar M}_{Col(\omega,<{\bar \kappa})} \
\{ X \in [H_{{\check {\bar \theta}}}^{V[{\dot G}]}]^\omega 
\colon X {\rm \ condenses \ remarkably \
} \} {\rm \ is \ stationary," }$$ inside ${\bar M}[G] \subset V$
we get some ${\bar \sigma} \colon 
H^{\bar M}_\rho \rightarrow H^{\bar M}_{\bar \theta}$ such that $\rho$
is a regular cardinal in ${\bar M}$ with $\rho < {\bar \kappa}$.

Now set $M = H_\rho^{\bar M}$, $N = H^{\bar M}_{\bar \theta}$,
$\sigma = {\bar \sigma}$, 
and $\pi =
{\bar \pi} \circ {\bar \sigma}$. Then $\pi$, $M$, $\sigma$, ${\bar \theta}$, and $N$
are as in \ref{def-remarkable}. As $\theta$ 
was arbitrary, this proves "$\Leftarrow$."

\bigskip
\hfill $\square$ (\ref{char-remarkable})

\begin{lemma}\label{rel-down-to-L}
Let $\kappa$ be remarkable. Then $L \models$ "$\kappa$ is remarkable."
\end{lemma}

{\sc Proof.} Let $\theta > \kappa$ be a regular cardinal in $L$. Let $G$ be
$Col(\omega,<\kappa)$-generic over $V$, and let
${\cal M} = (L_\theta[G];\in,{\vec R}) \in L[G]$ be any model of finite type. Let
${\cal N} = (H_{\theta^+}^{V[G]};\in,H_{\theta^+}^V,G,L_\theta[G],{\vec R})$.
As $\kappa$ is remarkable, in $V[G]$ we may pick some
$$\pi \colon 
(H_\beta^{V[G \cap H_\alpha]};\in,H_\beta^V,G \cap H_\alpha,
L_{\bar \theta}[G \cap L_\alpha],{\vec {\bar R}}) 
\rightarrow 
{\cal N}$$ where $\alpha = crit(\pi) < \beta < \kappa$ and $\beta$ is a regular
cardinal in $V$. Then $$\pi \upharpoonright L_{\bar \theta}[G \cap L_\alpha]
\colon (L_{\bar \theta}[G \cap L_\alpha],\in,{\vec {\bar R}})
\rightarrow {\cal M} {\rm , }$$ and ${\bar \theta}$ is a regular cardinal in $L$.
Because $L_{\bar \theta}[G \cap L_\alpha] \in L[G]$  and is countable there,
\ref{ex-of-emb2} and the existence of 
$\pi \upharpoonright L_{\bar \theta}[G \cap L_\alpha]$ yield
that inside $L[G]$ there are
predicates ${\vec {\bar S}}$ on   
$L_{\bar \theta}[G \cap L_\alpha]$ together with an elementary embedding
$$\sigma
\colon (L_{\bar \theta}[G \cap L_\alpha];\in,{\vec {\bar S}})
\rightarrow {\cal M}.$$ 
I.e., $(ran(\sigma);\in,{\vec R}
\upharpoonright ran(\sigma)) \prec {\cal M}$ where $ran(\sigma)
\in L[G]$. As $\theta$ and then ${\cal M}$ were
arbitrary we have shown that in $L$ does $\kappa$ satisfy the characterization of
remarkability from \ref{char-remarkable}.

\bigskip
\hfill $\square$ (\ref{rel-down-to-L})

\section{$L({\mathbb R})$ absoluteness.}

\begin{lemma}\label{lemma-on-ind-reals}
Let $\kappa$ be remarkable in $L$. Let $G$ be
$Col(\omega,<\kappa)$-generic over $L$. Let $P \in L[G]$ be a proper poset, and let
$H$ be $P$-generic over $L[G]$. Then for every real $x$ in $L[G][H]$ there is a poset
$Q_x \in L_\kappa$ such that $x$ is $Q_x$-generic over $L$.
\end{lemma}

{\sc Proof.} Let $\theta > \kappa$ be an $L$-cardinal such that ${\cal P}(P) \subset
L_\theta[G]$. Let $x \in {\mathbb R} \cap L[G][H]$, and let ${\dot x} \in
L_\theta[G]$ be such that ${\dot x}^H = x$.
Consider the structure ${\cal M} = 
(L_\theta[G];\in,P,{\dot x},H)$. Because $\kappa$ is remarkable in $L$ and
$P$ is proper we may pick an elementary
$$\pi \colon 
(L_\beta[G \cap L_\alpha];\in,{\bar P},{\bar {\dot x}},{\bar H})
\rightarrow {\cal M}$$ with the property that $G \cap L_\alpha$ is
$Col(\omega,<\alpha)$-generic over $L$ and $\beta$ is an $L[G \cap L_\alpha]$-cardinal.
By elementarity, ${\bar H}$ is ${\bar P}$-generic over $L_\beta[G \cap L_\alpha]$, and
hence over $L[G \cap L_\alpha]$, as ${\cal P}({\bar P}) \cap L[G \cap L_\alpha]
\subset L_\beta[G \cap L_\alpha]$. Moreover, by the definability of forcing,
we get that $n \in {\bar {\dot x}}^{\bar H}$ iff $\exists p \in {\bar H} \ p \forces
{\check n} \in {\bar {\dot x}}$ iff $\exists p \in H \ p \forces
{\check n} \in {\dot x}$ iff $n \in {\dot x}^H$ iff $n \in x$. So 
${\bar {\dot x}}^{\bar H} = x$, and we may set $Q_x = Col(\omega,<\alpha) \star {\dot
{\bar P}}$ where ${\dot {\bar P}}^{\bar H} = {\bar P}$. Notice that $Q_x \in L_\kappa$.

\bigskip
\hfill $\square$ (\ref{lemma-on-ind-reals}) 

\begin{lemma}\label{lemma-on-all-reals}
Let $\kappa$ be remarkable in $L$. Let $G$ be
$Col(\omega,<\kappa)$-generic over $L$. Let $P \in L[G]$ be a proper poset, and let
$H$ be $P$-generic over $L[G]$. Let $E$ be
$Col(\omega,(2^{\aleph_0})^{L[G][H]})$-generic over $L[G][H]$. Then in $L[G][H][E]$
there is some $G'$ being $Col(\omega,<\kappa)$-generic over $L$ such that $${\mathbb
R} \cap L[G'] = {\mathbb R} \cap L[G][H].$$ 
\end{lemma}

{\sc Proof.} Let $( e_i \colon i<\omega ) \in L[G][H][E]$ be such that
$\{ e_i \colon i<\omega \} = {\mathbb R} \cap L[G][H]$. By working inside $L[G][H][E]$
we may easily use \ref{lemma-on-ind-reals} to
construct $(\alpha_i,G_i \colon i<\omega)$ such that $\alpha_0 < \alpha_1 <
...$ and
for all $i<\omega$ we have that $G_i$ is $Col(\omega,<\alpha_i)$-generic over $L$,
$G_{i-1} \subset G_i$ (with the convention that $G_{-1} = \emptyset$), 
$G_i \in L[G][H]$, and 
$e_i \in L[G_i]$. Set $$G' = \bigcup_i \ G_i.$$ Because $Col(\omega,<\kappa)$ has the
$\kappa$-c.c., $G'$ is $Col(\omega,<\kappa)$-generic over $L$, and every real in
$L[G']$ is in $L[G_i]$ for some $i<\omega$. We get that ${\mathbb
R} \cap L[G'] = {\mathbb R} \cap L[G][H]$ as desired.

\bigskip
\hfill $\square$ (\ref{lemma-on-all-reals}) 

\begin{thm}\label{emb-thm-in-L[G]} {\em (Embedding theorem in $L[G]$)}
Let $\kappa$ be remarkable in $L$. Let $G$ be
$Col(\omega,<\kappa)$-generic over $L$, and write $V = L[G]$.  
Then in $V$ the $L({\mathbb R})$ embedding theorem holds for proper forcings.
\end{thm}

{\sc Proof.} Let $P \in V$ be a proper poset, and let $H$ be $P$-generic over $V$.
By \ref{lemma-on-all-reals} (in some further extension) there is $G'$
being $Col(\omega,<\kappa)$-generic over $L$ such that ${\mathbb
R} \cap L[G'] = {\mathbb R} \cap V[H]$. Let $\phi({\vec v},{\vec w})$ be a formula,
let ${\vec x} \in {\mathbb R} \cap V$, and let ${\vec \alpha} \in OR$. We then have
that $$L({\mathbb R}^V) \models \phi({\vec x},{\vec \alpha}) \Leftrightarrow$$
$$\forces^{L[{\vec x}]}_{Col(\omega,<\kappa)} \ L({\dot {\mathbb R}})
\models \phi({\check {\vec x}},{\check {\vec
\alpha}}) \Leftrightarrow$$ $$L({\mathbb R}^{V[H]}) \models 
\phi({\vec x},{\vec \alpha}).$$

\bigskip
\hfill $\square$ (\ref{emb-thm-in-L[G]})

\begin{thm}\label{anti-coding-in-L[G]} {\em (Anti-coding theorem in $L[G]$)}
Let $\kappa$ be remarkable in $L$. Let $G$ be
$Col(\omega,<\kappa)$-generic over $L$, and write $V = L[G]$.  
Then in $V$ the $L({\mathbb R})$ anti coding theorem holds for proper forcings.
\end{thm}

{\sc Proof.} Let $P \in V$ be a proper poset, and let $H$ be $P$-generic over $V$. 
By \ref{emb-thm-in-L[G]}, it suffices to show that each 
$A \in L({\mathbb R}^{V[H]}) \cap V$ is also in $L({\mathbb R}^V)$.
Fix such an $A$, and let $\Phi$ a formula, ${\vec \alpha} \in OR$, 
and $x \in {\mathbb R}^{V[H]}$ be such
that $$\gamma \in A \Leftrightarrow L({\mathbb R}^{V[H]}) \models 
\Phi({\vec \alpha},x,\gamma).$$ Let ${\dot x}^H = x$, and assume w.l.o.g. that
$$(\star) \ \ \ \ \ \emptyset \forces^V_P \ \gamma \in {\check A} \Leftrightarrow 
L({\dot {\mathbb R}}) \models 
\Phi({\check {\vec \alpha}},{\dot x},\gamma).$$ 
As in the proof of \ref{lemma-on-ind-reals},
we may pick an elementary  
$$\pi \colon 
(L_\beta[G \cap L_\alpha];\in,{\bar P},{\bar {\dot x}},{\bar H})
\rightarrow  
(L_\theta[G];\in,P,{\dot x},H)$$ such that $\beta$ is an 
$L[G \cap L_\alpha]$-cardinal. Because $L_\beta[G \cap L_\alpha]$ is countable in $V$
we may pick $h \in V$ being ${\bar P}$-generic over $L_\beta[G \cap L_\alpha]$. Of
course, $h$ will then also be ${\bar P}$-generic over $L[G \cap L_\alpha]$.
Because $P$ is proper we may and shall assume w.l.o.g. that (inside some further
forcing extension) for every $p \in {\bar P}$ there is $G^p$ being $P$-generic
over $V$ with $\pi(p) \in G^p$ 
and such that $\pi^{-1}{\rm " }G^p$ is ${\bar P}$-generic over 
$L_\beta[G \cap L_\alpha]$ (i.e., over $L[G \cap L_\alpha]$). Notice that ${\dot
x}^{G^p} = x$ for every $p \in {\bar P}$.
In order to prove \ref{anti-coding-in-L[G]} it now clearly suffices to verify the
following.

\bigskip
{\em Claim.} For all $\gamma \in OR$, $\gamma \in A \Leftrightarrow 
\forces^{L[G \cap L_\alpha][h]}_{Col(\omega,<\kappa)} \ L({\dot {\mathbb R}}) \models
\Phi({\check {\vec \alpha}},{\bar {\dot x}},{\check \gamma})$.

\bigskip
{\sc Proof.} We shall prove "$\Leftarrow$." The proof of "$\Rightarrow$" is almost
identical in that it starts from $\lnot \ \Phi$ instead of from $\Phi$, and gives
$\gamma \notin A$ instead of $\gamma \in A$. Suppose that 
$$\forces^{L[G \cap L_\alpha][h]}_{Col(\omega,<\kappa)} \ L({\dot {\mathbb R}}) \models
\Phi({\check {\vec \alpha}},{\bar {\dot x}}^h,{\check \gamma}).$$ This is itself 
forced by
some $p \in h$, and thus we also get, writing ${\bar G}^p = \pi^{-1}{\rm " }G^p$, that
$$\forces^{L[G \cap L_\alpha][{\bar G}^p]}_{Col(\omega,<\kappa)} 
\ L({\dot {\mathbb R}}) \models
\Phi({\check {\vec \alpha}},{\bar {\dot x}}^{{\bar G}^p},{\check \gamma}).$$
Because ${\bar G}^p = \pi^{-1}{\rm " }G^p \in L[G][G^p]$, in much the same way 
as in the proof of \ref{lemma-on-ind-reals} we can pick (inside some further forcing
extension) some $G'$ being $Col(\omega,<\kappa)$-generic over 
$L[G \cap L_\alpha][{\bar G}^p]$ such that 
$${\mathbb R} \cap L[G \cap L_\alpha][{\bar G}^p][G'] =
{\mathbb R} \cap L[G][G^p].$$ Hence $$L({\mathbb R}^{V[G^p]}) \models 
\Phi({\vec \alpha},{\bar {\dot x}}^{{\bar G}^p},\gamma).$$ But 
${\bar {\dot x}}^{{\bar G}^p} = {\dot x}^{G^p}$, so that there is some $q \in
G^p$ such that $$q \forces^V_P L({\dot {\mathbb R}}) \models \Phi({\check {\vec
\alpha}},{\dot x},{\check \gamma}).$$ Hence, by $(\star)$,
$q \forces^V_P {\check \gamma} \in {\check A}$, which implies that $\gamma \in A$.

\bigskip
\hfill $\square$ (Claim)

\hfill $\square$ (\ref{anti-coding-in-L[G]})

\bigskip
Here is an immediate corollary to \ref{emb-thm-in-L[G]} and
\ref{anti-coding-in-L[G]}, when combined with \ref{rel-down-to-L}.

\begin{cor}
Neither the conclusion of the $L({\mathbb R})$ embedding theorem 
for proper forcings nor the conclusion of the $L({\mathbb R})$
anti coding theorem for proper forcings implies $\Pi^1_1$-determinacy.
\end{cor}
\section{An equiconsistency.}

\begin{defn}
Let $A \subset OR$. We say that $A$ is {\em good} if $A \subset \omega_1$
and $L_{\omega_2}[A]
= H_{\omega_2}$.
\end{defn}

\begin{lemma}\label{good-A} 
If $0^\sharp$ does not exist then there is a proper $P \in V$ such that $$\forces_P \
{\rm "there \ is \ a \ good \ } A. {\rm " }$$
\end{lemma}

{\sc Proof.} This uses almost disjoint forcing in its
simplest form.
Fix $\delta$, a singular cardinal of uncountable cofinality and such
that $\delta^{\aleph_0} = \delta$ (for example, let $\delta$ be a
strong limit). By Jensen's Covering Lemma, we know that $\delta^{+L} =
\delta^+$. We may also 
assume w.l.o.g. that $2^\delta = \delta^+$, because
otherwise we may collapse $2^\delta$ onto $\delta^+$ by a
$\delta$-closed preliminary forcing. We may hence
pick $B \subset \delta^+$ with the property
that $H_{\delta^+} = L_{\delta^+}[B]$. 

Now let $G_1$ be $Col(\delta,\omega_1)$-generic over $V$. Notice
that the
forcing is $\omega$-closed. 
Set $V_1 = V[G_1]$. We have that
$\omega_2^{V_1} = \delta^+ = \delta^{+L}$. Let $C \subset \omega_1$
code $G_1$ (in the sense that $G_1 \in L_{\omega_2^{V_1}}[C]$).
Using the fact that $Col(\delta,\omega_1)$ has the 
$\delta^+$-c.c., it is easy to verify that in $V_1$, 
$H_{\omega_2} = L_{\omega_2}[B,C]$. Let $\omega_2$ denote $\omega_2^{V_1}$ from now on.

In $L$ we may pick $(A'_\xi \ \colon \ \xi < \delta^+)$, a
sequence of almost disjoint subsets of $\delta$. In
$L_{\omega_2}[C]$ we may pick a bijective $g \colon \omega_1
\rightarrow \delta$. Then if we let $\alpha \in A_\xi$ iff $g(\alpha)
\in A'_\xi$ for $\alpha < \omega_1$ and
$\xi < \delta^+$, we have that $(A_\xi \ 
\colon \ \xi < \delta^+)$ is a sequence of almost distoint subsets of
$\omega_1$. 

In $V_1$, we may pick $D \subset \omega_2$ with
$H_{\omega_2} = L_{\omega_2}[B,C] = L_{\omega_2}[D]$
(for example, the "join" of $A$ and $B$). We let $P_2$
be the forcing for coding $D$ by a subset of $\omega_1$, using the
almost disjoint sets $A_\xi$.

To be specific, $P_2$ consists of pairs $p = (l(p),r(p))$ where $l(p)
\colon \alpha \rightarrow 2$ for some $\alpha < \omega_1$ and $r(p)$
is a countable subset of $\omega_2$. We have $p = (l(p),r(p))
\leq_{P_2} q = (l(q),r(q))$ iff $l(p) \supset l(q)$, $r(p) \supset
r(q)$, and for all $\xi \in r(q)$, if $\xi \in D$ then $$\{ \beta
\in dom(l(p)) \setminus dom(l(q)) \ \colon \ l(p)(\beta) = 1 \} \cap
A_\xi \ = \ \emptyset.$$

By a $\Delta$-system argument, $P_2$ has the $\omega_2$-c.c. It is
clearly $\omega$-closed, so no cardinals are collapsed.
Moreover, if $G_2$ is $P_2$-generic over $V_1$, and if we set
$$A  \ = \ \bigcup_{ p \in G_2 } \ \{ \beta \in dom(l(p)) \ \colon \
l(p)(\beta) = 1 \} {\rm , }$$ then $A \subset \omega_1$
and we have that for all
$\xi < \omega_2$, $$\xi \in A_1 {\rm \ iff \ } 
Card(A \cap A_\xi) \leq
\aleph_0.$$

This means that $A_1$ is an element of any inner model containing
$(A_\xi \ \colon \ \xi < \omega_2)$ and $A$. (Of course, much more
holds.) An example of such a model is $L[D]$ in the sense explained
above. Set $V_2 = V_1[G_2]$. Moreover, 
because $P_2$ has the c.c.c., we get that in $V_2$,
$H_{\omega_2} = L_{\omega_2}[A]$.

Recall that all the forcings we have used to obtain $V_2$ were either $\omega$-closed
or had the c.c.c. Hence $V_2$ is a proper set-generic extension of $V$.

\bigskip
\hfill $\square$ (\ref{good-A})

\bigskip
It is easy to see that the conclusion of \ref{good-A} is actually equivalent with the
property that $V$ is not closed under $\sharp$'s.

\begin{defn}
Let $A \subset \omega_1$. By $\triangledown(A)$ we denote the assertion that  
$$\{ X \in [L_{\omega_2}[A]]^\omega \colon \exists \alpha < \beta \in
Card^{L[A \cap \alpha]} \ \exists \pi \ \pi \colon L_\beta[A \cap \alpha] \rightarrow 
X \prec
L_{\omega_2}[A] \}$$ is stationary in $[L_{\omega_2}[A]]^\omega$. 
\end{defn}

\begin{thm}\label{good-heart}
Suppose that $L({\mathbb R})$ is 
absolute under proper forcings.
Then $$\forall A \ (A {\rm \ good \ } \Rightarrow \triangledown(A))$$ 
holds in all proper set-generic extensions of $V$.
\end{thm}

{\sc Proof.} Let $\Psi$ denote the statement that the reals can be well-ordered in
$L({\mathbb R})$. By adding $\omega_1$ Cohen reals with finite support, which is
proper, one obtains an extension of $V$ in which 
$\Psi$ fails. Hence if $L({\mathbb R})$ is 
supposed to be
absolute under proper forcings, $\Psi$ has to fail in $V$ to begin with, and it has to
fail in every proper set-forcing extension of $V$.

Let us now fix a good $A$ such that 
$\triangledown(A)$ fails. We shall define a proper forcing
$P \in V$ such that $$\forces_P \ \Psi.$$ This will give a contradiction, and prove
\ref{good-heart} in $V$; of course, by replacing $V$ by a proper set-forcing extension
of itself, the very same argument will prove the full \ref{good-heart}.

The key observation here is that $\lnot \ \triangledown(A)$ implies that "reshaping" 
our good $A$ is proper. We let $P_1$ consist of functions $p \colon \alpha \rightarrow
2$ with $\alpha < \omega_1$ and such that for all $\xi \leq \alpha$ we have that $$L[A
\cap \xi,p \upharpoonright \xi] \models Card(\xi) \leq \aleph_0.$$ This is the
classical forcing for reshaping $A$ (cf. \cite{BeJeWe}). We need the following.

\bigskip
{\em Claim.} $P_1$ is proper. 

\bigskip
{\sc Proof.} Let us consider ${\cal M} =
(L_{\omega_2}[A];\in,A)$. Because $\triangledown(A)$
fails, there is a club $C \subset [H_{\omega_2}]^\omega$ 
such that for all $X \in C$, if 
$$\pi \colon
(L_\beta[A \cap \alpha];\in,A \cap \alpha) \cong (X;\in,A \cap X) \prec
{\cal M}$$ then 
$\beta$ is not a cardinal in $L[A \cap \alpha]$. Let us fix some such $X$. We have to
show that for any $p \in P_1 \cap X$ there is $q \leq_{P_1} p$ which is
$(P,X)$-generic.

For this we use an argument of \cite{ShSt}. Let $({\dot \alpha}_i \colon i<\omega)$
enumerate the ordinal names in $X$. We shall produce $q \leq_{P_1} p$ such that for
all $i<\omega$ we have that $q \forces {\dot \alpha} \in X$.
We may assume w.l.o.g. that $\alpha = \omega_1^{L[A \cap
\alpha]}$, as otherwise the 
task of constructing $q$ turns out to be an easier variant of what is to follow.
Now as $\beta$ has size 
$\alpha$ in $L[A \cap \alpha]$ 
we may pick a club $E \subset \alpha$ in
$L[A \cap \alpha]$ which grows faster than all clubs in $L_\beta[A \cap \alpha]$, 
i.e.,
whenever ${\bar E} \subset \alpha$ is a club in $L_\beta[A \cap \alpha]$ 
then $E \setminus
{\bar E}$ is bounded in $\alpha$. 

Inside $L[A \cap \alpha]$, 
we are now going to construct a sequence $(p_i \ \colon
\ i < \omega)$ of conditions below $p$
such that $p_{i+1} \leq_{P_1} p_i$ and
$p_{i+1} \forces {\dot \alpha}_i \in X$. We also want to maintain inductively
that $p_{i+1} \in L_\beta[A \cap \alpha]$. (Notice that $p \in 
L_\beta[A \cap \alpha]$
to begin with.) 
In the end we also want to have that setting $q =
\cup_{i < \omega} \ p_i$, we have that
$q \in P_1$, which of course is the the
non-trivial part. For this purpose, we also pick $({\bar \alpha}_i \colon i<\omega)$
cofinal in $\alpha$.

To commence, let $p_0 = p$. Now suppose that $p_i$ is given, $p_i \in
L_\beta[A \cap \alpha]$. Set
$\gamma = dom(p_i) < \alpha$. Work inside $L_\beta[A \cap \alpha]$ for a minute.
For all $\delta$ such that $\gamma \leq
\delta < \alpha$ we may pick some $p^\delta
\leq_{P_1} p_i$ such
that: $p^\delta \forces \pi^{-1}({\dot \alpha}_i) 
\in L_\beta[A \cap \alpha]$, $dom(p^\delta) \geq {\bar \alpha}_i$, 
$dom(p^\delta) > \delta$,
and for all limit ordinals
$\lambda$, $\gamma \leq \lambda \leq \delta$, $p^\delta(\lambda)
= 1$ iff $\lambda = \delta$.  
Then there is ${\bar E}$ club in $\alpha$ such that for any $\eta \in
{\bar E}$, $\delta < \eta \Rightarrow dom(p^\delta) <
\eta$. 

Now back in $L[A \cap \alpha]$, we may pick 
$\delta \in E$ such that $E \setminus {\bar E} \subset
\delta$. Set $p_{i+1} = p^\delta$, and let for future reference $\delta =
\delta_{i+1}$.
Of course $p_{i+1} \forces {\dot \alpha}_i \in X$.
We also have that $dom(p_{i+1}) < min \{ \epsilon \in E \ \colon \
\epsilon > \delta \}$, so that
for all limit ordinals $\lambda \in E
\cap (dom(p_{i+1}) \setminus dom(p_i))$ we have that 
$p_{i+1}(\lambda) = 1$ iff $\lambda = \delta_{i+1}$.

Now set $q = \cup_{i < \omega} \ p_i$. 
We are done if we can show that
$q$ is a condition. Well, it is easy to see that we have arranged that $dom(q) =
\alpha$, so that the only problem here is 
to show that $$L[A \cap \alpha,q]
\models Card(\kappa) \leq \aleph_0.$$ 
But by the construction
of the $p_i$'s we have that $$\{ \lambda \in E \cap (dom(q) \setminus
dom(p)) \ \colon \ \lambda {\rm \ is \ a \ limit \ ordinal \ and \ }
q(\lambda) = 1 \}$$ $$= \ 
\{ \delta_{i+1} \ \colon \ i < \omega \} {\rm ,
}$$ being a cofinal subset of $E$. But
$E$ is an element of $L[A \cap \alpha,q]$, so that
$\{ \delta_{i+1} \ \colon \ i < \omega \} \in L[A \cap \alpha,q]$
witnesses that $Card(\alpha) \leq \aleph_0$.

We have shown that $q \in P_1$ is $(P,X)$-generic, as desired.

\bigskip
\hfill $\square$ (Claim)

\bigskip
Now let $G$ be $P_1$-generic over $V$, and pick $D \subset \omega_1$ such that
$L_{\omega_2}[D] = L_{\omega_2}[A,G]$.
We may now "code down to a real" by using almost disjoint
forcing. By the fact that $D$ is "reshaped," there is a
(unique) sequence $(a_\beta \ \colon \ \beta 
< \omega_1)$ of 
subsets of $\omega$ such that for each $\beta 
< \omega_1$, $a_\beta$ is
the $L[D \cap \beta]$-least subset of $\omega$ being almost disjoint
from any $a_{\bar \beta}$ for ${\bar \beta} < \beta$.

We then let $P_2$ consist of all pairs $p = (l(p),r(p))$ where
$l(p) \colon n \rightarrow 2$ for some $n < \omega$ and $r(p)$ is a
finite subset of $\omega_1$. We let $p = (l(p),r(p)) \leq_{P_2} q =
(l(q),r(q))$ iff $l(p) \supset l(q)$, $r(p) \supset r(q)$, and for all
$\beta \in r(q)$, if $\beta \in D$ then $$\{ \gamma \in dom(r(p))
\setminus dom(r(q)) \ \colon \ r(p)(\gamma) = 1 \} \cap a_\beta \ = \
\emptyset.$$ 

By a $\Delta$-system argument, $P_1$ has the c.c.c..
Moreover, if $H$ is $P_2$-generic over $V[G]$, and if we set
$$a = \bigcup_{p \in H} \ \{ \gamma \in dom(l(p)) \ \colon \
l(p)(\gamma) = 1 \} {\rm , }$$ then we have that for $\gamma <
\omega_1$, $$\gamma \in D {\rm \ iff \ } Card(a \cap a_\gamma) <
\aleph_0.$$ 
Moreover, because $P_2$ has the c.c.c., we get that in $V[G][H]$ we have that
$H_{\omega_2} = L_{\omega_2}[a]$.

In particular, ${\mathbb R} \cap V[G][H] \subset L[a]$ which implies that in 
$V[G][H]$ there is a $\Delta^1_2(a)$-well-ordering of the reals. Thus, if we set $P =
P_1 \star {\dot P_2}$ then $P$ is proper and $$\forces_P \ \Psi.$$

\bigskip
\hfill $\square$ (\ref{good-heart})   

\begin{lemma}\label{abs-heart}
Suppose that $L({\mathbb R})$ is 
absolute under proper forcings.
Then $\omega_1$ is remarkable in $L$.
\end{lemma}

{\sc Proof.} By \ref{ind-rem}, we may assume that $0^\sharp$ does not exist.
Let $\theta > \kappa$ be an $L$-cardinal. Using \ref{good-A} we may easily find a
proper set-forcing extension of $V$ in which there is a good $B$ and in which $\theta
< \omega_2$ (just primarily 
force with $Col(\omega_1,\theta)$, which is $\omega$-closed).
By finally forcing with $Col(\omega,<\omega_1)$ (which has the c.c.c.) we get a 
proper set-forcing extension of $V$ in which we may pick a good 
$A$ such that $A_{odd} = \{ 2 \delta + 1 \in A \colon \delta < \omega_1 \}$
essentially is $Col(\omega,<\omega_1)$-generic over $L$, and $\theta < \omega_2$.
By \ref{good-heart} we know that in that extension,
$$\{ X \in [L_{\omega_2}[A]]^\omega \colon \exists \alpha < \beta \in
Card^{L[A \cap \alpha]} \ \exists \pi \ \pi \colon L_\beta[A \cap \alpha] \rightarrow 
X \prec
L_{\omega_2}[A] \}$$ is stationary in $[L_{\omega_2}[A]]^\omega$. We may now argue 
exactly as
in the proof of \ref{rel-down-to-L} to see that this implies that $\omega_1$ has to be
remarkable in $L$.

\bigskip
\hfill $\square$ (\ref{abs-heart})   

\begin{cor}\label{equicon}
The following are equiconsistent.

(1) $L({\mathbb R})$ is 
absolute under proper forcings.

(2) The $L({\mathbb R})$ embedding theorem holds for proper forcings.

(3) The $L({\mathbb R})$ anti coding theorem holds for proper forcings.

(4) There is a remarkable cardinal.
\end{cor}

{\sc Proof.} ${\sc Con}(1) \Rightarrow {\sc Con}(4)$ is \ref{abs-heart}. 
${\sc Con}(4) \Rightarrow {\sc Con}(2)$ and
${\sc Con}(4) \Rightarrow {\sc Con}(3)$ are \ref{emb-thm-in-L[G]} and
\ref{anti-coding-in-L[G]}. ${\sc Con}(3) \Rightarrow {\sc Con}(4)$
follows from the proofs of \ref{good-A} and \ref{good-heart}.
${\sc Con}(2) \Rightarrow {\sc Con}(1)$ is trivial.

\bigskip
\hfill $\square$ (\ref{equicon})   

\section{A derived model theorem.}

We have shown in \ref{emb-thm-in-L[G]} that there is a model of $L({\mathbb R})$
absoluteness under proper forcings which is of the form $L[G]$ where $G$ is
$Col(\omega,<\kappa)$-generic over $L$ for some inaccessible $\kappa$ in $L$.
We are now going to show that -- under some genericity assumption -- every 
model of $L({\mathbb R})$
absoluteness under proper forcings is of this form.

\begin{defn}\label{genericity}
We let $(\natural)$ denote the assertion that every real is set-generic over $L$, i.e.,
that for every $x \in {\mathbb R}$ there is some poset $P \in L$ and some $G \in V$ being
$P$-generic over $L$ such that $x \in L[G]$.
\end{defn}

\begin{thm}\label{DMT} {\em (Derived model theorem)}
Assume that $(\natural)$ holds and that $L({\mathbb R})$ is absolute under proper
forcings. Then (in some set-generic extension of $V$) there is $G$ being
$Col(\omega,<\omega_1^V)$-generic over $L$ such that $L({\mathbb R}^V) =
L({\mathbb R}^{L[G]})$.
\end{thm}

{\sc Proof.} By \ref{good-A} and \ref{good-heart} there is $V[H]$, a proper set-generic
extension of $V$, in which there is a good $A$, and $\triangledown(A)$ holds.
By $(\natural)$, for every $x \in {\mathbb R}^V$ we may pick
a poset $P_x \in L$ and some $K_x \in V$ being $P_x$-generic over $L$ 
such that $x \in L[K_x]$. Let $\theta_x$ be such that $P_x \in H_{\theta_x}$.
By primarily forcing with $Col(\omega_1,sup_{x \in
{\mathbb R}}(Card(P)))$ 
we may assume w.l.o.g. that any $P_x$ is hereditarily smaller than $\omega_2$ in
$V[H]$, i.e.,
$P_x \in L_{\omega_2}[A]$ for every $x \in {\mathbb R}^V$.
 
Now fix $x \in {\mathbb R}^V$, and set ${\cal M} = (L_{\omega_2}[A];\in,A,P_x,K_x,{\dot
x})$ where ${\dot x}^{K_x} = x$. Using $\triangledown(A)$
there is some $$\pi \colon (L_\beta[A \cap \alpha];\in,A \cap \alpha,{\bar P_x},{\bar
K_x},{\bar {\dot x}}) \rightarrow {\cal M}$$ such that $\beta$ is a cardinal in 
$L[A \cap \alpha]$, and hence so in $L$. We get that $x = ({\bar {\dot x}})^{\bar K_x} 
\in L[{\bar
K_x}]$ where ${\bar K_x}$ is
${\bar P_x}$-generic over $L$, and ${\bar P_x}$ is countable.   
Notice that $\pi$ only exists in $V[H]$. However, by \ref{ex-of-emb2}
we may then also find, inside $V$, some 
$$\sigma \colon (L_\beta;\in,{\tilde P_x},{\tilde K_x},{\tilde {\dot x}}) 
\rightarrow {\cal M} {\rm , }$$ so that $x = ({\tilde {\dot x}})^{\tilde K_x} 
\in L[{\tilde
K_x}]$ where ${\tilde K_x}$ is
${\tilde P_x}$-generic over $L$, and ${\tilde P_x}$ is countable.

But now, as in the proof of \ref{lemma-on-all-reals}, in a
$Col(\omega,(2^{\aleph_0})^V)$-generic extension of $V$ we may construct $G$ being
$Col(\omega,<\omega_1^V)$-generic over $L$ such that $L({\mathbb R}^V) =
L({\mathbb R}^{L[G]})$.

\bigskip
\hfill $\square$ (\ref{DMT})

\begin{cor}
Assume that $(\natural)$ holds and that $L({\mathbb R})$ is absolute under proper
forcings. Then every set of reals in $L({\mathbb R})$ is Lebesgue measurable and has
the property of Baire.
\end{cor}

\end{document}